\documentclass[letterpaper]{amsart}

\usepackage[T1]{fontenc}
\usepackage{amssymb}
\usepackage{amsmath}
\usepackage{amsthm}
\usepackage[shortlabels]{enumitem}
\usepackage{mathtools}
\usepackage{xcolor}
\usepackage{listings}
\usepackage[obeyspaces]{url}
\usepackage{hyperref}

\numberwithin{equation}{section}
\numberwithin{figure}{section}

\theoremstyle{plain}
\newtheorem{theorem}{Theorem}[section]

\theoremstyle{definition}

\newtheorem{remark}[theorem]{Remark}

\newcommand{\R}{\mathbb{R}}
\newcommand{\C}{\mathbb{C}}

\newcommand{\Sphere}{\widehat{\mathbb{C}}}
\newcommand{\CodeSphere}{\smash[t]{\Sphere}}
\newcommand{\cF}{\mathcal{F}}
\newcommand{\loc}{\mathrm{loc}}

\newcommand{\Mod}{\operatorname{Mod}}

\let\lean\nolinkurl
\newcommand{\leanfile}[1]{\href{https://github.com/will1491/no-wandering-domains/blob/1a39faf/NoWanderingDomains/#1}{\lean{#1}}}
\newcommand{\leanfolder}[1]{\href{https://github.com/will1491/no-wandering-domains/tree/1a39faf/NoWanderingDomains/#1}{\lean{#1}}}

\definecolor{LeanCodeBackground}{HTML}{F5F5F2}
\definecolor{LeanCodeBorder}{HTML}{E2E2DC}
\definecolor{LeanKeyword}{HTML}{B4232B}
\definecolor{LeanType}{HTML}{087F72}
\definecolor{LeanBuiltin}{HTML}{2457C5}
\definecolor{LeanComment}{HTML}{6B7280}
\definecolor{LeanString}{HTML}{287A3D}

\lstdefinelanguage{Lean}{
  sensitive=true,
  alsoletter={_'},
  morekeywords=[1]{
    abbrev,axiom,class,def,deriving,do,else,end,example,extends,
    fun,if,import,in,inductive,infix,infixl,infixr,instance,let,
    macro,match,namespace,noncomputable,opaque,open,private,
    protected,section,structure,then,theorem,universe,variable,
    where,with
  },
  morekeywords=[2]{
    Type,Prop,Sort,Set,Finset,Nat,Int,Real,Complex,
    BeltramiCoeff,Measurable
  },
  morekeywords=[3]{
    apply,atTop,by,cases,constructor,exact,fderiv,intro,intros,
    infer_instance,nhds,obtain,refine,rfl,rw,simp,simpa,tsupport,
    volume
  },
  morecomment=[l]{--},
  morecomment=[s]{/-}{-/},
  morestring=[b]"
}

\lstdefinestyle{LeanCodeStyle}{
  language=Lean,
  basicstyle=\small\ttfamily,
  keywordstyle=[1]\color{LeanKeyword}\bfseries,
  keywordstyle=[2]\color{LeanType},
  keywordstyle=[3]\color{LeanBuiltin},
  commentstyle=\color{LeanComment}\itshape,
  stringstyle=\color{LeanString},
  backgroundcolor=\color{LeanCodeBackground},
  rulecolor=\color{LeanCodeBorder},
  frame=single,
  framerule=0.4pt,
  framesep=7pt,
  aboveskip=0.8\baselineskip,
  belowskip=0.8\baselineskip,
  xleftmargin=0pt,
  xrightmargin=0pt,
  columns=fullflexible,
  keepspaces=true,
  showstringspaces=false,
  breaklines=true,
  breakatwhitespace=false,
  mathescape=true,
  literate={:=}{{{\color{LeanBuiltin}:=}}}2
           {=>}{{{\color{LeanBuiltin}=>}}}2,
  tabsize=2
}

\lstnewenvironment{leancode}
  {\lstset{style=LeanCodeStyle}}
  {}

\title[No Wandering Domains in Lean]{Formalization of Sullivan's No Wandering Domains Theorem in Lean}
\author{Ziang Li}
\author{Yusheng Luo}
\date{September 10, 2026}


\begin{document}

\begin{abstract}
We report on a Lean 4 formalization of Sullivan's No Wandering Domains theorem: every Fatou component of a rational self-map of the Riemann sphere of degree at least two is eventually periodic. The project formalizes the relevant background in normal families, the Montel--Carath\'eodory theorem, Julia and Fatou sets, local Sobolev regularity and Wirtinger derivatives, the Cauchy and Beurling transforms, the equivalence of analytic and geometric quasiconformality, and the measurable Riemann mapping theorem. This paper describes the translation of this mathematics into Lean, the proof architecture, the reusable components, and the autoformalization workflow.
\end{abstract}

\maketitle

\section{Introduction}

Let $f\colon \Sphere\to\Sphere$ be a rational map of degree $d\geq 2$. The Fatou set $\Omega(f)$ is the locus on which the sequence of iterates $\{f^n\}_{n\geq 0}$ is normal, and the Julia set $J(f)$ is its complement. A connected component $U$ of $\Omega(f)$ is \emph{wandering} if the components containing $f^n(U)$ are pairwise disjoint. Sullivan's theorem says that this never occurs: every Fatou component is eventually periodic \cite{Sullivan}. The result is central in modern complex dynamics because it reduces the study of all Fatou components to periodic components and their preimages.

The mathematical statement is short, but its formal proof requires a long chain of analytic and dynamical infrastructure. Normal-family theory underlies the definitions of the Fatou and Julia sets and the proofs of their basic properties. The deformation argument uses measurable Beltrami differentials, weak derivatives, and solutions of the $\bar\partial$-equation. The broader project also proves the equivalence of analytic and geometric quasiconformality and the measurable Riemann mapping theorem, although these final results lie outside the dependency chain of the No Wandering Domains theorem. Formalizing these topics required substantial background in real analysis, measure theory, topology, and complex analysis. To the best of our knowledge, many of the specialized results needed here were not available in Mathlib or other existing Lean projects when the development began.

The repository develops this infrastructure on top of Mathlib, RMT4, and Carleson \cite{Carleson,RMT4,Mathlib}. The source code is available at \url{https://github.com/will1491/no-wandering-domains}. At commit \lean{1a39faf}, using Lean and Mathlib \lean{v4.33.0}, it contains 171 Lean files and approximately 157,000 lines under \lean{NoWanderingDomains/}. The project source contains no uses of \lean{sorry} or \lean{admit} and introduces no project-specific \lean{axiom} declarations.

The formalized proof follows the development in McMullen's course notes \cite{McMullen}. If $U$ were wandering, one could construct a space of $f$-invariant Beltrami differentials of arbitrarily large dimension. We choose dimension $2d+2$ for simplicity. Solving the $\bar\partial$-equation defines a linear map from this space to the $(2d+1)$-dimensional space of infinitesimal holomorphic deformations of $f$. Rank--nullity therefore gives a nonzero invariant Beltrami differential whose induced deformation vanishes. The remainder of the proof shows that such a differential must itself vanish, yielding a contradiction.

This paper has two goals beyond reporting our formalization. First, we aim to make the Lean infrastructure we built more reusable in future projects. By describing in detail the definitions, tools, and theorems we developed along the way to the No Wandering Domains theorem in Sections~\ref{sec:statements} and~\ref{sec:components}, we hope to give readers a clearer understanding of how to use this infrastructure in their own projects. Second, to the best of our knowledge, this project is among the early formalizations in complex dynamics. As formal verification evolves rapidly, we want to illustrate what autoformalization can currently achieve with appropriate human guidance in this field and, more generally, in formalizing difficult mathematical theorems.

\section{Mathematical statements and Lean formulations}
\label{sec:statements}

We present the Lean formulations of the No Wandering Domains theorem and the measurable Riemann mapping theorem, along with the definitions required for stating them. Links to the corresponding source files are provided for reference.

\subsection{No wandering domains}

\begin{theorem}[Sullivan's No Wandering Domains theorem]
If $f\colon\Sphere\to\Sphere$ is rational of degree at least two, then no connected component of $\Omega(f)$ is wandering. Equivalently, every Fatou component is eventually periodic.
\end{theorem}

The Lean formulation in \leanfile{Dynamics/NoWanderingDomains.lean} is:

\begin{leancode}
theorem sullivan_no_wandering_domains
    {f : $\CodeSphere \to \CodeSphere$}
    (hf : IsRational f)
    (hd : 2 $\leq$ degreeOfRational f) :
    $\forall$ U : Set $\CodeSphere$,
      IsFatouComponent f U $\to$ $\neg$ IsWandering f U
\end{leancode}

The Lean statement uses \lean{IsFatouComponent f U} and \lean{IsWandering f U}. A Fatou component is a connected component of the Fatou set $\Omega(f)$, as defined in \leanfile{Dynamics/FatouComponents/Def.lean}. It is wandering if the Fatou components of its forward orbits $f^n(U)$ are pairwise distinct, as formalized in \leanfile{Dynamics/FatouComponents/Periodic.lean}.

For $f\colon\Sphere\to\Sphere$, the Fatou set $\Omega(f)$ is the locus where the family of iterates is normal, and the Julia set $J(f)$ is its complement. In \leanfile{Dynamics/JuliaFatou/Def.lean}:

\begin{leancode}
def FatouSet (f : $\CodeSphere \to \CodeSphere$) : Set $\CodeSphere$ :=
  {z | IsNormalAt (Set.range fun n : $\N$ => f^[n]) z}

def JuliaSet (f : $\CodeSphere \to \CodeSphere$) : Set $\CodeSphere$ :=
  (FatouSet f)${}^{\mathrm c}$
\end{leancode}

The condition \lean{IsNormalAt} is defined through normal families. A family $\cF\subseteq X\to Y$ is normal on $U\subseteq X$ if every sequence in $\cF$ has a subsequence converging locally uniformly on $U$. It is normal at a point if it is normal on some neighborhood. For reuse, \leanfile{NormalFamilies/Basic.lean} defines these notions for any topological space $X$ and uniform space $Y$:

\begin{leancode}
def IsNormal {X Y : Type*}
    [TopologicalSpace X] [UniformSpace Y]
    (F : Set (X $\to$ Y)) (U : Set X) : Prop :=
  $\forall$ seq : $\N$ $\to$ F,
    $\exists$ phi : $\N$ $\to$ $\N$, StrictMono phi $\land$
    $\exists$ g : X $\to$ Y,
      TendstoLocallyUniformlyOn
        (fun n => (seq (phi n) : X $\to$ Y)) g atTop U

def IsNormalAt {X Y : Type*}
    [TopologicalSpace X] [UniformSpace Y]
    (F : Set (X $\to$ Y)) (z : X) : Prop :=
  $\exists$ U $\in$ nhds z, IsNormal F U
\end{leancode}

The normal-family infrastructure also includes the Montel--Carath\'eodory theorem: a family of sphere-holomorphic maps omitting three fixed distinct points is normal for the spherical metric \cite{MontelCaratheodory}. Its Lean statement, \lean{montel_caratheodory_sphere}, is in \leanfile{NormalFamilies/StrongMontel/SphereMontel.lean}.

\subsection{The measurable Riemann mapping theorem}

\begin{theorem}[Measurable Riemann mapping theorem]
For every Beltrami coefficient $\mu$, there exists an orientation-preserving quasiconformal homeomorphism $f\colon\C\to\C$ satisfying $\bar\partial f=\mu\partial f$ almost everywhere. There is a unique such solution satisfying $f(0)=0$ and $f(1)=1$ \cite{AhlforsBers}.
\end{theorem}

The existence statement is formalized in \leanfile{QC/MRMT/Existence.lean} as follows:

\begin{leancode}
theorem mrmt_exists (b : BeltramiCoeff) :
  $\exists$ f : $\C \to \C$, IsQCAnalytic f b
\end{leancode}

The normalized uniqueness statement is formalized in \leanfile{QC/MRMT/Uniqueness.lean}:

\begin{leancode}
theorem mrmt_unique_normalized (b : BeltramiCoeff) :
  $\exists!$ f : $\C \to \C$,
    IsQCAnalytic f b $\land$ f 0 = 0 $\land$ f 1 = 1
\end{leancode}

The Lean statements use \lean{BeltramiCoeff} for the coefficient and \lean{IsQCAnalytic} for the solution.

A Beltrami coefficient is a measurable function $\mu\colon\C\to\C$ with $\|\mu\|_{\infty}<1$. In \leanfile{QC/Defs/BeltramiCoeff.lean}:

\begin{leancode}
structure BeltramiCoeff where
  $\mu$ : $\C \to \C$
  measurable : Measurable $\mu$
  bound : eLpNormEssSup $\mu$ volume < 1
\end{leancode}

The analytic definition of quasiconformality requires $f$ to be an orientation-preserving homeomorphism in $W^{1,2}_{\loc}(\C)$ satisfying the Beltrami equation $\bar\partial f=\mu\partial f$ almost everywhere. Here, $W^{1,2}_{\loc}(\C)$ means that $f$ and its first weak derivatives are locally square-integrable, encoded by \lean{MemW12loc} in \leanfile{Analysis/Sobolev/WeakDeriv.lean}. The symbols $\partial f$ and $\bar\partial f$ denote the Wirtinger derivatives, implemented as \lean{dz} and \lean{dzbar} in \leanfile{Analysis/Sobolev/Wirtinger.lean}. Orientation preservation is expressed by positivity of the Jacobian almost everywhere. In \leanfile{QC/Defs/Analytic.lean}:

\begin{leancode}
def OrientationPreservingHomeo (f : $\C \to \C$) : Prop :=
  IsHomeomorph f $\land$
    $\forall^{\mathrm{a.e.}}$ z, 0 < (fderiv $\R$ f z).det

def IsQCAnalytic (f : $\C \to \C$) (b : BeltramiCoeff) : Prop :=
  OrientationPreservingHomeo f $\land$ MemW12loc f $\land$
    $\forall^{\mathrm{a.e.}}$ z, dzbar f z = b.$\mu$ z * dz f z
\end{leancode}

The project also formalizes the geometric definition of quasiconformality and proves its equivalence with the analytic definition. This allows future formalizations to use whichever definition is more natural and transfer results between the two.

The geometric definition begins with the orientation condition \lean{SensePreserving} defined in \leanfile{QC/Defs/SensePreserving.lean}: $f$ is a homeomorphism that sends sufficiently small positively oriented circles around almost every $z_0$ to curves of winding number $+1$ about $f(z_0)$.

A quadrilateral $Q$ is represented by a continuous map $\R^2\to\C$ that is injective on the closed unit square. Write $D_Q$ for the image of the square and $L_Q$ and $R_Q$ for the images of its left and right sides. The connecting family $\Gamma(Q)$ consists of all absolutely continuous curves $\gamma\colon[0,1]\to D_Q$ with $\gamma(0)\in L_Q$ and $\gamma(1)\in R_Q$. Its conformal modulus is
\[
  \Mod(\Gamma(Q))=
  \inf_{\rho\ \mathrm{admissible}}
  \int_{\C}\rho(z)^2\,dA(z),
\]
where admissibility means that $\int_\gamma\rho\,ds\geq1$ for every $\gamma\in\Gamma(Q)$.

For a homeomorphism $f$, \lean{Q.imageCurveFamily f}, denoted by $\Gamma(fQ)$, consists of all absolutely continuous curves in $f(D_Q)$ joining $f(L_Q)$ to $f(R_Q)$. The geometric definition requires
\[
  \Mod(\Gamma(fQ))\leq K\,\Mod(\Gamma(Q))
\]
for every quadrilateral $Q$. Geometric quasiconformality is encoded in \leanfile{QC/Defs/Geometric.lean} as follows:

\noindent\begin{minipage}{\linewidth}
\begin{leancode}
def IsQCGeometric (f : $\C \to \C$) (K : $\R$) : Prop :=
  1 $\leq$ K $\land$ SensePreserving f $\land$
    $\forall$ Q : Quadrilateral,
      curveModulus (Q.imageCurveFamily f)
        $\leq$ ENNReal.ofReal K * Q.modulus
\end{leancode}
\end{minipage}

\begin{remark}
To the best of our knowledge, no Lean formalization of the Jordan curve theorem was available at the time of this project. We therefore define a quadrilateral by a parametrization of the unit square, which gives its region, distinguished sides, and connecting curve family directly, without needing to recover a domain from its boundary.
\end{remark}

The equivalence between the two definitions of quasiconformal mapping is formalized in \leanfile{QC/Equivalence.lean}:

\begin{leancode}
theorem qc_analytic_iff_geometric
    {f : $\C \to \C$} {K : $\R$} (hK : 1 $\leq$ K) :
  ($\exists$ b : BeltramiCoeff,
      b.normInf $\leq$ (K - 1) / (K + 1) $\land$
      IsQCAnalytic f b)
    $\leftrightarrow$ IsQCGeometric f K
\end{leancode}

\section{Outline of the proof}
\label{sec:proofs}

The formalized proof of Sullivan's No Wandering Domains theorem proceeds in six steps. Let $f\colon\Sphere\to\Sphere$ be a rational map of degree $d\geq2$, and suppose for contradiction that $U$ is a wandering Fatou component.

\begin{enumerate}[label=\textbf{Step \arabic*.},leftmargin=*]
\item \textbf{Prepare the wandering component.} After replacing $U$ by a component farther along its forward orbit, one may assume that every iterate is injective on $U$ and that the orbit of $U$ contains no critical points.

\item \textbf{Construct invariant Beltrami differentials.} Choose a disk $D=B(a,\rho)$ with $\overline D\subset U$ and, for $0\leq k\leq2d+1$, define
\[
  \nu_k(z)=(k+1)\,\overline{(z-a)}^k\,\mathbf{1}_D(z).
\]
Their span is
\[
  E=\operatorname{span}_{\C}\{\nu_0,\ldots,\nu_{2d+1}\}, \qquad \dim_{\C}E=2d+2.
\]
As in McMullen's proof of Theorem~5.33 \cite{McMullen}, each $\nu\in E$ can be propagated along the grand orbit of $U$ to an $f$-invariant Beltrami differential $\mu_\nu$ on $\Sphere$. The differential $\mu_\nu$ agrees with $\nu$ on $U$ and vanishes outside the grand orbit of $D$. The disjointness of the orbit components and the injectivity of the iterates on $U$ make this propagation well defined, and the map $\nu\mapsto\mu_\nu$ is linear.

\item \textbf{Associate an infinitesimal deformation.} Solve $\bar\partial v_\nu=\mu_\nu$ using the linear $\bar\partial$-solution operator (infinitesimal MRMT). The invariance of $\mu_\nu$ implies that
\[
  \eta_\nu=f'v_\nu-v_\nu\circ f
\]
is a holomorphic section of $f^*T\Sphere$. Thus, the construction defines a linear map
\[
  E\longrightarrow H^0(\Sphere,f^*T\Sphere),
  \qquad \nu\longmapsto\eta_\nu.
\]

\item \textbf{Apply the dimension count.} Since $T\Sphere\simeq\mathcal O(2)$ and $f$ has degree $d$, we have $f^*T\Sphere\simeq\mathcal O(2d)$ and therefore
\[
  \dim_{\C}H^0(\Sphere,f^*T\Sphere)=2d+1.
\]
The domain $E$ has dimension $2d+2$, so rank--nullity gives a nonzero $\nu_*\in E$ for which $\eta_{\nu_*}=0$. Equivalently,
\[
  v_{\nu_*}(f(z))=f'(z)v_{\nu_*}(z).
\]

\item \textbf{Use repelling periodic points.} If $p$ is a repelling periodic point of period $n$, then iteration of the preceding equation gives
\[
  v_{\nu_*}(p)=(f^n)'(p)v_{\nu_*}(p).
\]
Since $|(f^n)'(p)|>1$, it follows that $v_{\nu_*}(p)=0$. Repelling periodic points are dense in $J(f)$, so continuity gives $v_{\nu_*}=0$ on $J(f)$ and hence on $\partial U$.

\item \textbf{Obtain the contradiction.} Write $\nu_*=\sum_{k=0}^{2d+1}c_k\nu_k$. For each $\nu_k$, an explicit solution of $\bar\partial\Phi_k=\nu_k$ is given by $\Phi_k(z)=\overline{(z-a)}^{k+1}$ on $D$ and $\Phi_k(z)=\rho^{2(k+1)}(z-a)^{-(k+1)}$ outside $D$. Define $h=\sum_{k=0}^{2d+1}c_k\Phi_k-v_{\nu_*}$.
Since $U$ is wandering, $\mu_{\nu_*}$ agrees with $\nu_*$ on $U$. 
Thus, the function $h$ is holomorphic on $U$. Outside $D$, the sum of the $\Phi_k$ is the rational function
\[
  R(z)=\sum_{k=0}^{2d+1}c_k\rho^{2(k+1)}(z-a)^{-(k+1)},
\]
so $h=R-v_{\nu_*}$ on $U\setminus\overline D$. Since $v_{\nu_*}$ vanishes on $\partial U$, the functions $h$ and $R$ have the same boundary values. The boundary-uniqueness theorem therefore forces every $c_k$ to vanish. This contradicts $\nu_*\neq0$, and therefore $f$ has no wandering Fatou component.
\end{enumerate}

\section{Reusable components}
\label{sec:components}

In addition to the headline results highlighted in Section~\ref{sec:statements}, the project develops substantial analytic and dynamical machinery that may be useful for future formalization projects and mathematical research. This section presents several of these reusable components.

\subsection{The Cauchy transform}

For a compactly supported function $\omega$, the Cauchy transform is defined by
\[
  P\omega(z)=-\frac{1}{\pi}\int_{\C}
    \frac{\omega(\zeta)}{\zeta-z}\,dA(\zeta).
\]
The transform is used twice in the project: to construct solutions of the Beltrami equation in the proof of the measurable Riemann mapping theorem, and to solve $\bar\partial v=\mu$ in the proof of No Wandering Domains.

The folders \leanfolder{Analysis/SingularIntegral} and \leanfolder{QC/MRMT/NeumannSeries} contain the smooth and $L^p$ theories of this transform. The formalization proves $P(\bar\partial\omega)=\omega$ and $\bar\partial(P\omega)=\omega$ for $C^1$ functions of compact support. For compactly supported $h\in L^p(\C)$ with $p>2$, it proves that $Ph$ is H\"older continuous with exponent $1-2/p$, tends to zero at infinity, and satisfies $\bar\partial(Ph)=h$ and $\partial(Ph)=Th$ weakly.

\subsection{The Beurling transform}

The Beurling transform is the principal-value singular integral
\[
  T\mu(z)=-\frac{1}{\pi}\operatorname{p.v.}
    \int_{\C}\frac{\mu(\zeta)}{(z-\zeta)^2}\,dA(\zeta).
\]
The transform is used to construct solutions of the Beltrami equation in the proof of the measurable Riemann mapping theorem. Since $T=\partial\circ P$, writing $f(z)=z+Ph(z)$ reduces the equation $\bar\partial f=\mu\,\partial f$ to
\[
  h=\mu(1+Th).
\]

The folder \leanfolder{Analysis/SingularIntegral/Beurling} develops the transform from the truncated Calder\'on--Zygmund operators of the Carleson project \cite{Carleson}. It proves $T=\partial\circ P$ and establishes
\begin{align*}
  \|T\mu\|_2&=\|\mu\|_2,\\
  \|T\mu\|_p&\leq C\|\mu\|_p \qquad (1<p<\infty).
\end{align*}

\subsection{Weak compactness and Jacobian limits}

Given a measurable Beltrami coefficient $\mu$, the project approximates $\mu$ by smooth functions, solves the corresponding smooth Beltrami equations, and passes to a limit. The folder \leanfolder{Analysis/WeakLimits} develops the weak $L^2$ compactness and Jacobian convergence results needed for this passage. In \leanfolder{QC/Calculus}, these results show that a locally uniform limit of geometric $K$-quasiconformal homeomorphisms remains $K$-quasiconformal when the limit is again a homeomorphism.

\section{Code development and autoformalization}
\label{sec:development}

The code was developed through autoformalization with Claude Code, using Fable 5 and Opus 4.8. The mathematical architecture was designed by the authors: the target theorem was decomposed into normal-family, hyperbolic, Sobolev, singular-integral, quasiconformal, MRMT, and deformation layers, and interfaces between those layers were chosen before the final assembly. The language models were used to turn these blueprints into Lean declarations and proofs, to search the existing library, and to revise the code in response to errors reported by Lean.

Every generated proof must ultimately be accepted by Lean's kernel. If the generated code contains an invalid proof step, Lean rejects it, allowing the error to be identified and corrected. However, kernel checking does not decide whether a definition captures the intended mathematics, whether a theorem has been stated at a useful level of generality, or whether an important hypothesis has been encoded in an unnatural or unnecessarily restrictive way. The authors carefully checked and verified these mathematical choices throughout the development.

After the formalization of the No Wandering Domains theorem was completed, the code underwent manual cleanup to improve file organization, code conciseness and cleanliness, and docstring readability. This included reorganizing files by dependency, removing duplicated local arguments, shortening proofs after stable helper lemmas emerged, reducing imports, strengthening docstrings, standardizing names, and separating public endpoints from implementation details.

\end{document}